\documentstyle[12pt]{ article}
\begin{document}

\newtheorem{theorem}{Theorem}[section]
\newtheorem{prop}[theorem]{Theorem}
\newenvironment{mtheorem}{\begin{prop}\rm}{\end{prop}}
\newtheorem{coroll}[theorem]{Corollary}
\newenvironment{mcorollary}{\begin{coroll}\rm}{\end{coroll}}
\newtheorem{lemm}[theorem]{Lemma}
\newenvironment{mlemma}{\begin{lemm}\rm}{\end{lemm}}
\newtheorem{anoprop}[theorem]{Proposition}
\newenvironment{mproposition}{\begin{anoprop}\rm}{\end{anoprop}}
\newenvironment{mproof}{\begin{trivlist}\item[]{\em
Proof: }}{\hfill$\Box$\end{trivlist}}
\newenvironment{mdefinition}{\begin{trivlist}\item[]
{\em Definition.}}{\hfill$\Box$\end{trivlist}}
\newtheorem{eg}{\rm\sl \uppercase{Example}}[section]
\newenvironment{example}{\begin{eg}\rm}{\hfill$\Box$\end{eg}}

\newcommand{\ca}{{\cal A}}
\newcommand{\cb}{{\cal B}}
\newcommand{\cc}{{\cal C}}
\newcommand{\cd}{{\cal D}}
\newcommand{\cf}{{\cal F}}
\newcommand{\cg}{{\cal G}}
\newcommand{\ch}{{\cal H}}
\newcommand{\ck}{{\cal K}}
\newcommand{\cl}{{\cal L}}
\newcommand{\cn}{{\cal N}}
\newcommand{\cm}{{\cal M}}
\newcommand{\co}{{\cal O}}
\newcommand{\cs}{{\cal S}}
\newcommand{\ct}{{\cal T}}

\newcommand{\csa}{{$C^*$-algebra}}
\newcommand{\nsa}{{non-self-adjoint}}
\newcommand{\cssa}{{$C^*$-subalgebra}}

\def \IR{\hbox{{\rm I}\kern-.2em\hbox{{\rm R}}}}
\def \iR{\hbox{{\sevenrm I\kern-.2em\hbox{\sevenrm R}}}}
\def \IN{\hbox{{\rm I}\kern-.2em\hbox{\rm N}}}
\def \IC{\hbox{{\rm I}\kern-.6em\hbox{\bf C}}}
\def \IQ{\hbox{{\rm I}\kern-.6em\hbox{\bf Q}}}
\def \ZZ{\hbox{{\rm Z}\kern-.4em\hbox{\rm Z}}}

\newcommand{\ga}{{\alpha}}
\newcommand{\gb}{{\beta}}
\newcommand{\gd}{{\delta}}
\newcommand{\gA}{{\Alpha}}
\newcommand{\gB}{{\Beta}}
\newcommand{\gG}{{\Gamma}}
\newcommand{\gD}{{\Delta}}
\newcommand{\gs}{{\sigma}}
\newcommand{\gS}{{\Sigma}}

\newcommand{\ot}{\otimes}
\newcommand{\op}{\oplus}
\newcommand{\pr}{\prime}

\newcommand{\bmi}{\bf $i$}
\newcommand{\bmj}{\bf $j$}
\renewcommand{\varphi}{\phi}
\parindent0cm

\begin{center}
{\Large\bf
On the C*-envelope of approximately
finite-dimensional operator algebras}
\end{center}
\vspace{.3in}

\begin{center}
{\Large  C. Laurie*}\\
\it Mathematics Department\\
University of Alabama\\
USA
\vspace{.3in}
\rm
\end{center}

\begin{center}
{\Large S. C. Power}\\
\it Department of Mathematics and Statistics\\
Lancaster University\\
England LA1 4YF
\rm
\vspace{1cm}
\end{center}

\sf
\begin{abstract}
The C*-envelope of the limit algebra (or limit space) of a contractive regular
system of digraph algebras (or digraph spaces)
is shown to be an approximately finite C*-algebra and the direct system
for the C*-envelope is determined explicitly.
\end{abstract}
\rm
\vspace{1in}
* Partially supported by an NSF-EPSCoR travel grant.

\newpage

A number of recent studies of non-self-adjoint operator algebras have
been concerned with the Banach algebra direct limits of direct systems

\begin{center}
\setlength{\unitlength}{0.0125in}%
\begin{picture}(225,22)(125,755)
\thicklines
\put(240,760){\vector( 1, 0){ 40}}
\put(210,755){\makebox(0,0)[lb]{\raisebox{0pt}[0pt][0pt]{$A_2$}}}
\put(125,755){\makebox(0,0)[lb]{\raisebox{0pt}[0pt][0pt]{$A_1$}}}
\put(150,760){\vector( 1, 0){ 40}}
\put(310,755){\makebox(0,0)[lb]{\raisebox{0pt}[0pt][0pt]{\twlit . . . . }}}
\put(250,765){\makebox(0,0)[lb]{\raisebox{0pt}[0pt][0pt]{$\phi_2$}}}
\put(350,755){\makebox(0,0)[lb]{\raisebox{0pt}[0pt][0pt]{$ A$}}}
\put(160,765){\makebox(0,0)[lb]{\raisebox{0pt}[0pt][0pt]{$\phi_1$}}}
\end{picture}
\end{center}

in which the building block algebras are finite-dimensional digraph
algebras. The most accessible of these arise when the algebra
homomorphisms $\phi_k$ are regular
in the sense that (partial) matrix unit systems
can be chosen for $A_1, A_2, \dots$ so that each $\phi_k$ maps matrix
units to sums of matrix units. Under these conditions the maps
$\varphi_k$ need not be star-extendible or isometric
 and the C*-envelope of the limit
algebra, in the sense of Hamana, need not be an approximately
finite-dimensional C*-algebra. (See \cite{scp-simplicial}.)

In what follows we show that if the homomorphisms are regular
and contractive then the C*-envelope \ $C^*_{env} (A)$ \ of the limit
algebra \ $A$\  is an AF C*-algebra. Furthermore we identify explicitly a
direct system of finite-dimensional C*-algebras for $C^*_{env} (A)$. In
particular, this generalizes the result for
the triangular compression limit
algebras $\displaystyle{\lim_\to (T_{n_k}, \varphi_k)}$ considered
by Hopenwaser and Laurie \cite{hop-lau} and answers the problem posed
there.
 The relationship between $A$
and $C^*_{env} (A)$ is rather subtle ; even after a natural
telescoping of the given direct system to an essentially
isometric system (Lemma 9) it need not be the case that
$C^*_{env} (A)$ is an isometric direct limit of the algebras $C^* (\varphi_k
(A_k))$.

In the first section we consider only finite-dimensional matters. In
particular the  contractive regular morphisms are characterised as the
regular bimodule maps of compression type. In the second section we
obtain the main result and discuss a variety of examples.

\newpage

{\Large \bf 1. Regular contractive morphisms}
\vspace{.3in}

Let $G$ be a finite directed graph with no multiple edges and with
vertices labelled $1, 2, \dots, n$. Let $\{e_{ij} : 1 \leq i, j \leq
n\}$ be a matrix unit system for the full complex matrix algebra $M_n$.
Define $A(G)$ to be the linear span of those matrix units $e_{ij}$ for
which $(i,j)$ is an edge of $G$. If $G$ is a reflexive digraph then we refer to
$A(G)$ as a {\it digraph space}.
 If $G$ is reflexive and transitive (as a binary relation) then we
refer to $A(G)$ as a {\it digraph algebra}.
 In intrinsic terms a digraph algebra
(alias finite-dimensional CSL algebra/poset algebra/ incidence algebra)
is a subalgebra of $M_n$ containing a maximal abelian self-adjoint
subalgebra. Usually we consider the digraph algebras that are associated
with the standard matrix unit system for $M_n$.
\vspace{.3in}

{\bf Definition 1.} \ Let $A(G), A(H)$ be digraph spaces with associated
(standard) matrix unit systems. Then a linear map $\phi : A(G) \rightarrow
A(H)$
 is
said to be a  {\it regular bimodule map}, with respect to the given
matrix unit systems, if $\varphi$ maps matrix units to orthogonal sums of
matrix units and $\varphi$ is a bimodule map with respect to the
standard diagonal subalgebras.
That is, if $C(G)$ and $C(H)$ are these diagonals,
 then $\varphi (C(G)) \subseteq
C(H)$ and
\[
\phi (c_1 a c_2) = \phi (c_1) \phi (a) \phi (c_2)
\]
for all $c_1, c_2$ in $C(G)$ and $a$ in $A(G)$.
\vspace{.3in}

More generally a linear map between digraph spaces is said to be a
{\it regular bimodule map} if matrix unit systems can be chosen so that the
map is of the form above.

The terminology above should be compared with the following
more general
terminology (which is not needed in this paper).
A map $\alpha$ between digraph algebras is
said to be regular if partial matrix unit systems can be chosen  so that
for the associated diagonals, say $C$ and $D$ respectively,
$\alpha$
is a $C-D$ bimodule map and is regular in the sense that $\alpha$
maps the normaliser of $C$ into the normaliser of $D$. In the case of
contractive maps it can be shown that
 this notion coincides with that which is given
in the definition above.
That this notion is more general can be seen by considering the
(Schur) automorphisms of the 4-cycle digraph algebra which leave the diagonal
invariant.

If $\varphi$ is as in Definition 1 then the image $\varphi (e_{ij})$ of
each matrix unit $e_{ij}$ in $A(G)$ is, by assumption, a partial
isometry which is an orthogonal sum of matrix units. Note also that the
initial and final projections of $\varphi (e_{ij})$ are dominated
(perhaps properly) by the
diagonal projections $\varphi (e_{jj})$ and $\varphi (e_{ii})$
respectively.
\vspace{.3in}

{\bf Definition 2.}\  A diagonal projection $Q$ in $A(G)$ is said to be
{\it $A(G)$-irreducible} if
$C^* (Q A(G) Q) = QM_nQ = B(Q \ \IC^n)$.
\vspace{.3in}

Plainly $Q$ is $A(G)$-irreducible if the graph of the digraph algebra
$QA(G)Q$ on $Q\ \IC^n$ is connected as an
undirected graph. Also,
if $Q$ is a diagonal projection of $A(G)$, then, by considering
connected components of the graph of $QA(G)Q$, we can write $Q=\sum Q_i$, with
$Q_i$ diagonal and $A(G)$-irreducible, such that
\[
QAQ = \sum_i \oplus Q_i A Q_i
\]
\vspace{.3in}

{\bf Definition 3.} The map $\phi:  A(G) \rightarrow A(H)$ is
a {\it (standard) elementary
compression type map} associated with diagonal projections $(Q^G, P^H)$
if

(i) $Q^G$ and $P^H$ are diagonal projections in $A(G)$ and $A(H)$,
respectively, with rank $Q^G$ = rank $P^H = r,$ and

(ii) $Q^G$ is $A(G)$-irreducible, and

(iii) $\phi$ is a linear map of the form $\phi = \beta \circ \alpha$ where
$\alpha: A(G) \rightarrow M_r$ is compression by $Q^G$
and $\beta: M_r \rightarrow C^* (A(H))$ is a C*-algebra injection,
with $\beta(I) = P^H$, which is a linear
extension of a correspondence
$e_{ij} \rightarrow f_{{n_i}{n_j}}$
of standard matrix units.

A (standard)  elementary compression type map is simply the
regular bimodule map associated with an identification
of a connected full subgraph of $G$ with an
isomorphic subgraph of $H$. The minimal diagonal subprojections of
$P^H$ are associated with the vertices of the subgraph of $H$.
Note that $\phi$ depends
not only on the pair $(Q^G,P^H)$ and the matrix unit system
but also on the
the particular identification $\beta$ chosen for the pair. This dependence
is often suppressed in the subsequent discussions.

\vspace{1em}
A map $\phi: A(G) \rightarrow A(H)$ is a {\it (standard) compression type
map} if it is a
direct sum of elementary compression type maps.

More generally, $\phi$ is a {\it compression type map}
 if matrix unit systems can
be chosen so that $\phi$ is of compression
type with respect to these systems.
Compression type maps feature in the linear topological aspects
of limit algebras discussed in \cite{scp-ban}.

\vspace{.3in}
{\bf Lemma 4.}\  \it Let $\phi$ be the Schur product
projection map which is defined on
the digraph space of an m-cycle by deleting an entry corresponding
to a proper edge. That is, consider the map

\[
A \in \left[ \begin{array}{lcc} \begin{array}{ccc}
 $*$ & $*$ & \\ & $*$ & $*$
\end{array}  & & \\ & \ddots & \\
\begin{array}{cc}  & \\ $*$ & \end{array}
 & & \begin{array}{cc} $*$ & $*$\\ & $*$
\end{array}
\end{array}
\right]
\rightarrow \phi(A) \in
 \left[ \begin{array}{lcc} \begin{array}{ccc} * & * & \\ & * & *
\end{array}  & & \\ & \ddots & \\
\begin{array}{cc}  & \\ 0 & \end{array}
 & & \begin{array}{cc} * & *\\ & *
\end{array}
\end{array} \right]
\]

where the entry $a_{m1}$ of $A$ is replaced by  zero.
Then the norm of the linear map $\phi$ dominates
the quotient

\[
\frac{\cos (\pi /(2m+1))}{\cos (\pi /2m)}.
\]

In particular, this sparse triangular truncation map is not contractive.
\rm
\vspace{.3in}

\begin{mproof}
Consider the following matrix $A$ where, as usual, unspecified entries
are zero :

\[
A =
 \left[ \begin{array}{lcc} \begin{array}{ccc} 1 & 1 & \\ & 1 & 1
\end{array}  & & \\ & \ddots & \\
 \begin{array}{ll}  &  \\\!\!\!\!\!-1& \end{array} & &
 \begin{array}{cc} 1&1\\ &1 \end{array}
\end{array}
 \right]
\ \ \ \ \ (m \times m).
\]

Let $D = \hbox{ diag} \{1, w, w^2, \dots, w^{m-1}\}$ where $w =
\hbox{exp} (\pi i/m)$.
Since \mbox{$- \bar{w}^{m-1} =w$} we have

\[
D^* AD = I + wS = \left[
\begin{array}{lcc} \begin{array}{ccc} 1 & w & \\ & 1 & w \end{array} & & \\
& \ddots & \\
\begin{array}{ll} &  \\ \! w &  \end{array} & &
\begin{array}{cc} 1 & w \\ & 1 \end{array}
 \end{array}
 \right].
\]

where $S$ is the cyclic backward shift.
Spectral theory then yields that

\[
||I+wS||=|1+w|=2 \cos
(\frac{\pi}{2m}).
\]

 On the other hand the truncate of $A$ has norm equal
to $2 \cos (\frac{\pi}{2 m+1})$.
For details of this see, for example, Example 1.2.5 in \cite{goo-del-jon}.
\end{mproof}
\vspace{.3in}

{\bf Theorem 5.} \it Let $\phi : A(G) \rightarrow A(H)$ be a regular
bimodule map between digraph spaces.
Then $\phi$ is contractive if and only if $\phi$ is of compression type.
\rm
\vspace{.3in}

\begin{mproof}
Let $\phi$ be a regular bimodule map
with respect to the matrix unit systems $\{e_{ij}\}$
for $A(G)$ and $\{f_{kl}\}$ for $A(H)$.
Choose $f_{kk}$ in $A(H)$ such that $f_{kk}$ is a summand of $\phi
(e_{jj})$ for some $j$. Let $p_1 = f_{kk}$. Let ${\cal O} = \{p_1, p_2,
\dots, p_s\}$ be the set of distinct minimal diagonal projections in $A(H)$
which is  the orbit of $p_1$ under $\{ \phi (e_{ij}) : (i,j) \in
E(G)\}$.
This is the smallest set of projections which contains $p_1$ and is such
that if $p$ is in the set and $\varphi (e_{ij}) p \neq 0$ (or $p \varphi
(e_{ij}) \neq 0)$ then $\varphi (e_{ij}) p (\varphi (e_{ij}))^*$ (or
$\varphi (e_{ij})^* p \varphi e_{ij})$ also belongs to the set.
Each $p_i$ is a summand of $\phi (e_{{k_i}{k_i}})$ for some $k_i$.

We claim
that $k_i \neq k_j$ for $i \neq j$.

Suppose $k_i = k_j$ for some $i \neq j$. Then there is a subset $\{q_1,
\dots, q_{l +1}\}$ of the orbit with $q_1 = p_i$ and $q_{l+1} =
p_j$ which corresponds to a cycle in $G$ consisting of vertices
$\{n_1, \dots, n_l\}$ (with $n_1 = k_i$) and (directed) edges $E_m$
connecting $n_m$ and $n_{m+1} \hbox{ for }m=1, \dots ,l-1$, and edge $E_l$,
connecting $n_l$ and $n_1$ such that $\phi (E_m)$ maps $q_m$ to $q_{m+1}$
(or $q_{m+1}$ to $q_m$ depending on the direction of  $E_m$),
for $ m=1, \dots, l$.
Here we write $E_m$ for  the matrix
unit determined by the edge $E_m$.

By relabelling the vertices of $G$ appropriately we can assume
that $\{ n_1, \dots, n_l \}$ is $ \{1, \dots, l\}$ and that the edge $E_l$
runs from $1$ to $l$. For example, the cycle might look like

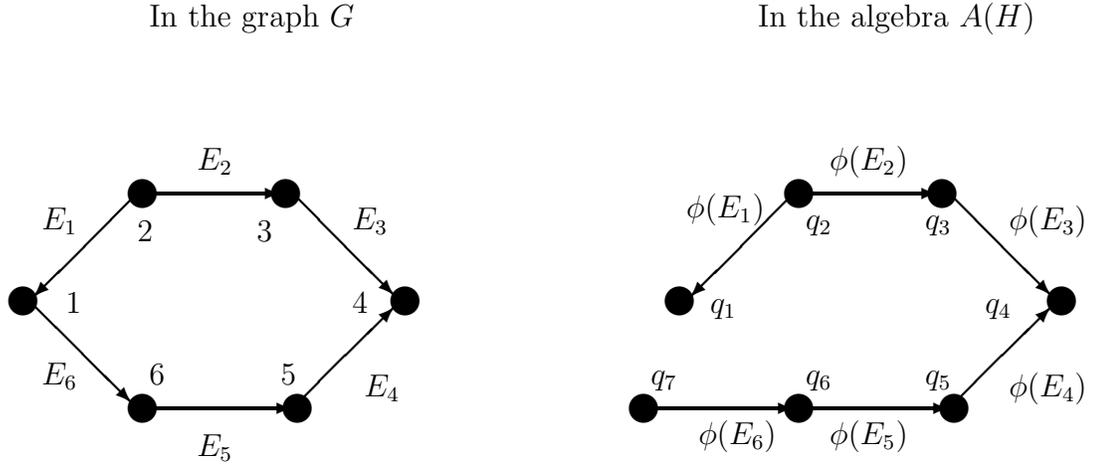
\begin{center}
\setlength{\unitlength}{0.0125in}%
\begin{picture}(449,192)(116,545)
\thicklines
\put(237,565){\circle*{12}}
\put(172,565){\circle*{12}}
\put(122,610){\circle*{12}}
\put(172,655){\circle*{12}}
\put(232,655){\circle*{12}}
\put(170,655){\vector( 1, 0){ 60}}
\put(170,655){\vector(-1,-1){ 43.500}}
\put(235,565){\vector( 1, 1){ 42.500}}
\put(235,655){\vector( 1,-1){ 42.500}}
\put(175,565){\vector( 1, 0){ 60}}
\put(125,610){\vector( 1,-1){ 42.500}}
\put(140,605){\makebox(0,0)[lb]{\raisebox{0pt}[0pt][0pt]{\twlrm 1}}}
\put(220,635){\makebox(0,0)[lb]{\raisebox{0pt}[0pt][0pt]{\twlrm 3}}}
\put(260,605){\makebox(0,0)[lb]{\raisebox{0pt}[0pt][0pt]{\twlrm 4}}}
\put(230,575){\makebox(0,0)[lb]{\raisebox{0pt}[0pt][0pt]{\twlrm 5}}}
\put(175,575){\makebox(0,0)[lb]{\raisebox{0pt}[0pt][0pt]{\twlrm 6}}}
\put(170,635){\makebox(0,0)[lb]{\raisebox{0pt}[0pt][0pt]{\twlrm 2}}}
\put(130,640){\makebox(0,0)[lb]{\raisebox{0pt}[0pt][0pt]{$ E_1$}}}
\put(195,665){\makebox(0,0)[lb]{\raisebox{0pt}[0pt][0pt]{$ E_2$}}}
\put(260,640){\makebox(0,0)[lb]{\raisebox{0pt}[0pt][0pt]{$ E_3$}}}
\put(265,570){\makebox(0,0)[lb]{\raisebox{0pt}[0pt][0pt]{$ E_4$}}}
\put(195,545){\makebox(0,0)[lb]{\raisebox{0pt}[0pt][0pt]{$ E_5$}}}
\put(130,575){\makebox(0,0)[lb]{\raisebox{0pt}[0pt][0pt]{$ E_6$}}}
\put(382,565){\circle*{12}}
\put(557,610){\circle*{12}}
\put(282,610){\circle*{12}}
\put(512,565){\circle*{12}}
\put(175,725){\makebox(0,0)[lb]{\raisebox{0pt}[0pt][0pt]{\twlrm In the graph
 $G$}}}
\put(447,565){\circle*{12}}
\put(397,610){\circle*{12}}
\put(447,655){\circle*{12}}
\put(507,655){\circle*{12}}
\put(445,655){\vector( 1, 0){ 60}}
\put(445,655){\vector(-1,-1){ 43.500}}
\put(510,565){\vector( 1, 1){ 42.500}}
\put(510,655){\vector( 1,-1){ 42.500}}
\put(450,565){\vector( 1, 0){ 60}}
\put(385,565){\vector( 1, 0){ 60}}
\put(450,640){\makebox(0,0)[lb]{\raisebox{0pt}[0pt][0pt]{$ q_2$}}}
\put(500,640){\makebox(0,0)[lb]{\raisebox{0pt}[0pt][0pt]{$ q_3$}}}
\put(525,605){\makebox(0,0)[lb]{\raisebox{0pt}[0pt][0pt]{$ q_4$}}}
\put(500,575){\makebox(0,0)[lb]{\raisebox{0pt}[0pt][0pt]{$ q_5$}}}
\put(450,575){\makebox(0,0)[lb]{\raisebox{0pt}[0pt][0pt]{$ q_6$}}}
\put(460,665){\makebox(0,0)[lb]{\raisebox{0pt}[0pt][0pt]{$ \phi(E_2$)}}}
\put(400,645){\makebox(0,0)[lb]{\raisebox{0pt}[0pt][0pt]{$ \phi(E_1$)}}}
\put(535,640){\makebox(0,0)[lb]{\raisebox{0pt}[0pt][0pt]{$ \phi(E_3$)}}}
\put(535,570){\makebox(0,0)[lb]{\raisebox{0pt}[0pt][0pt]{$ \phi(E_4$)}}}
\put(460,550){\makebox(0,0)[lb]{\raisebox{0pt}[0pt][0pt]{$ \phi(E_5$)}}}
\put(405,550){\makebox(0,0)[lb]{\raisebox{0pt}[0pt][0pt]{$ \phi(E_6$)}}}
\put(410,605){\makebox(0,0)[lb]{\raisebox{0pt}[0pt][0pt]{$ q_1$}}}
\put(385,575){\makebox(0,0)[lb]{\raisebox{0pt}[0pt][0pt]{$ q_7$}}}
\put(430,725){\makebox(0,0)[lb]{\raisebox{0pt}[0pt][0pt]{\twlrm In the algebra
 $A(H)$}}}
\end{picture}
\end{center}

Let $$a = \sum^{l-1}_{m=1} E_m + \sum_{t \in T} e_{tt} -
e_{l1}$$ where

\[
T = \{m: E_m \hbox{ and } E_{m-1} \hbox{ run in the same direction, where
} E_0 = E_l\}.
\]

In the  illustrated example,
\[ a = e_{12} + e_{32} + e_{43} + e_{45} + e_{56} +
e_{11} + e_{33} + e_{55} + e_{66} - e_{61}.
\]
Then, after deleting rows of zeros and columns of zeros,
$a$ has a submatrix of the form
\[
A =
 \left[ \begin{array}{lcc} \begin{array}{ccc} 1 & 1 & \\ & 1 & 1
\end{array}  & & \\ & \ddots & \\
 \begin{array}{ll} &  \\\!\!\!\!\! -1 &  \end{array} & &
 \begin{array}{cc} 1&1\\ &1 \end{array}
\end{array} \right].\]

Let $p = q_1 + \dots + q_\ell$. Since $\varphi (e_{\ell 1})$ maps
$q_{\ell +1}$ into $q_\ell, \ p \varphi (e_{\ell 1}) p=0$ and hence $p
\varphi (a) p$ has the same norm as the associated matrix

\[
 \left[ \begin{array}{lcc} \begin{array}{ccc} 1 & 1 & \\ & 1 & 1
\end{array}  & & \\ & \ddots & \\
 \begin{array}{cc} &  \\ 0 & \end{array} & &
 \begin{array}{cc} 1&1\\ &1 \end{array}
\end{array} \right].\]

By Lemma 3, $\parallel  a \mid
\mid < \parallel p \varphi (a) p \parallel $  which contradicts
the hypothesis that
$\varphi$ is contractive. Thus the claim is proven. That is, for each $p_i$
in the orbit  $\cal O$ with $p_i$ a summand in $\varphi (e_{k_i
k_i}),$ we have that $k_i \neq k_j$ if $p_i \neq p_j$.

Let $P = \sum p_i$, the sum taken over all $p_i$ in $\co$. Then $P$ is a
diagonal projection in $A(H)$. Let $Q = \sum e_{k_i k_i}$, the sum
taken over all $k_i$ such that $p_i$ is in $\cal O$. By the claim, $Q$
is a diagonal projection in $A(G)$ with rank $Q$ = rank $P$. Note that
the graph of $QA(G)Q$ is connected so $Q$ is $A(G)$-irreducible. Note
also, from the definition of $P$,
 that $(I - P) \varphi (A(G)) P=0$ and $P \varphi (A(G))(I
-P)=0$.

Thus we can write $\phi = \gamma \oplus \phi^\prime$ where
$\gamma =P\phi$ and $\phi^\prime = P^\bot\phi.$

There are two possibilities for $\gamma$. Either $\gamma$
is an elementary compression
type map (of multiplicity one) or $\gamma$ fails to be injective. In the latter
case
there is a non-self-adjoint
matrix unit, $e_{ij}$ say, which is mapped to zero by $\gamma$ (and $\phi$) and
which corresponds to an edge $E$ of $G$, which is not one of the edges
$E_1, \dots \ E_l$, but which nevertheless has its two vertices in common with
two of the vertices of $E_1, \dots \ E_l$.
This means that $E$, together with some of the edges
$E_1, \dots \ E_l$, form a cycle. The argument above applies and once again
we obtain the contradiction $\|\gamma\| > 1$.

Induction completes the proof.
\end{mproof}
\vspace{.3in}

{\bf Corollary 6.} \it If $\varphi : A(G) \rightarrow A(H)$ is a contractive
regular bimodule map, then $\varphi$ is completely contractive.
\rm
\vspace{.3in}

\begin{mproof}
For $\varphi$ of compression type, $\varphi^{(n)} : A(G)
\otimes M_n \rightarrow A(H) \otimes M_n$ is also of compression type.
\end{mproof}

We now examine further properties of compression type maps.
Let  $B(\cm)$
denote all bounded operators on the Hilbert space $\cm$ and write $\cong$ for
C*-algebra isomorphism.
\vspace{.3in}

{\bf Proposition 7.} \ \it Let $A(G) \subseteq B(\IC^n)$ and
$A(H) \subseteq B(\IC^m)$ be digraph algebras and let
$\varphi : A(G) \rightarrow A(H)$ be a
compression type
map.
Thus
$\varphi = \sum \op \gamma_i$ where each $\gamma_i$ is an
elementary compression type map associated with the pair $(Q^G_i, P_i^H)$
(and an implicit identification) with
each $Q_i^G$ being $A(G)$-irreducible.
Then

\[
C^* (\varphi (A(G))) \cong \sum_{k \in K} \oplus B(Q^G_k \ \IC^n)
\cong \sum_{k \in K} \oplus B(P_k^H\  \IC^m)
\]

where $\{Q^G_k : k
\in K\}$ includes exactly one copy of each distinct $Q^G_i$.
\vspace{.3in}
\rm

\begin{mproof}
We have  $\gamma_i = \beta_i \circ \alpha_i$ where $\alpha_i : A(G)
\rightarrow B(Q^G_i \ \IC^n)$
is the compression map $a \to Q_i^GaQ_i^G$
and $\beta_i : B(Q^G_i \ \IC^n) \rightarrow B
(P_i^H \ \IC^m)$ is a
(standard) C*-algebra isomorphism. Define
\[
\alpha : A(G) \rightarrow \sum_i \oplus
B(Q_i^G \ \IC^n)
\]
and
\[
\beta : \sum \oplus B(Q_i^G \ \IC^n) \rightarrow \sum \oplus B(P_i^H
\ \IC^m)
\]
by $\alpha(a) = \sum \oplus \ \alpha_i (a)$ and
$\beta (\sum \oplus b_i) = \sum \oplus \beta_i (b_i)$.
Thus $C^* (\varphi(A(G)))$ is
isomorphic to $C^* (\alpha (A(G)))$
which in turn is isomorphic to $C^* (\hat{\alpha}
(A(G)))$ where $\hat{\alpha}(a) = \sum_{i \in K} \op \alpha_i (a)$.

Thus $C^*(\phi(A(G))$ is identified with a C*-subalgebra, $E$ say,
of
\[
F = \sum_{k \in K} \oplus B(Q^G_k \ \IC^n).
\]
The compression of $E$ to each summand of $F$ is equal to the
summand, so it remains to show only that no summand
of $E$ appears with multiplicity
in the summands of $F$. This is elementary.
For example if $Q_k^G$ and $Q_l^G$ are distinct,
with $k, l \in K,$ and \ rank($Q_k^G$) = rank($Q_l^G$)\ then
consider a self-adjoint matrix unit $e$ of $A(G)$ with $Q_k^Ge = e$
and $Q_l^Ge = 0$. Then $\phi(e)$ differs in the summands for
$Q_k^G$ and $Q_l^G$ being nonzero in one
summand and zero in the other, as desired.
\end{mproof}

\vspace{.3in}

{\bf Proposition 8.} \it Let $\gamma : A(G) \rightarrow A(H)$ be an
elementary compression type map  with
projection pair $(Q^G, P^H)$ and let
$\eta : A(H) \rightarrow A(F)$ be an elementary compression type map
with projection pair $(Q^H, P^F)$.
Then \\ \mbox{$\eta \circ \gamma : A(G) \rightarrow
A(F)$} is of compression type with $\eta \circ \gamma = \sum \op \delta_i$
where $\delta_i$ is an elementary compression type map  with
projection pair $(q_i^G, p_i^F)$ where
the $q^G_i$ are orthogonal subprojections of $Q^G$ and
the  $p_i^F$ are orthogonal subprojections of $P^F$.
\vspace{.3in}
\rm

\begin{mproof}
Suppose $P^H = \sum_{k \in K} f_{kk}$ so that $P^H Q^H = \sum_{j \in J}
f_{jj}$ for some subset $J$ of $ K$. Each $f_{jj}, j \in J$, corresponds via
$\gamma$ with a minimal diagonal projection $e_{{k_j}{k_j}}$ which is
a summand of $Q^G$. Let $q^G = \sum_{j \in J} e_{{k_j}{k_j}}$ and
write $q^G = \sum q_i^G$ as a direct sum of $A(G)$-irreducible
projections. Each $q_i^G$ corresponds under $\gamma$ to a diagonal
subprojection, $p_i^H$, of $P^H Q^H$, which in turn corresponds under
$\eta$ to a diagonal subprojection $p_i^F$ of $P^F$. Thus $\eta \circ
\gamma = \sum \delta_i$ where $\delta_i$ is an elementary compression
type map associated with $(q_i^G, p_i^F)$.
\end{mproof}
\vspace{.3in}

{\large \bf 2. Regular direct systems of digraph spaces}
\rm

We turn our attention now to systems

\begin{center}
\setlength{\unitlength}{0.0125in}%
\begin{picture}(225,22)(125,755)
\thicklines
\put(240,760){\vector( 1, 0){ 40}}
\put(210,755){\makebox(0,0)[lb]{\raisebox{0pt}[0pt][0pt]{$A_2$}}}
\put(125,755){\makebox(0,0)[lb]{\raisebox{0pt}[0pt][0pt]{$A_1$}}}
\put(150,760){\vector( 1, 0){ 40}}
\put(310,755){\makebox(0,0)[lb]{\raisebox{0pt}[0pt][0pt]{\twlit . . . . }}}
\put(250,765){\makebox(0,0)[lb]{\raisebox{0pt}[0pt][0pt]{$\phi_2$}}}
\put(160,765){\makebox(0,0)[lb]{\raisebox{0pt}[0pt][0pt]{$\phi_1$}}}
\end{picture}
\end{center}

where
each $A_i$ is a digraph space $A (G_i)$
for some digraph $G_i$ and each $\varphi_i$ is a
contractive regular bimodule map.
We refer to such a system as a {\it contractive regular direct
system of digraph spaces}.

Define
\[
A^0_\infty = \{(a_k) : a_k
\in A_k, \varphi_k (a_k) = a_{k+1} \mbox{ for all large } k\}
\]
and let $A_\infty$ be
the set of equivalence classes of eventually equal sequences. Let
$\varphi_{k, \infty}$ denote the natural map from $A_k$ into $A_\infty$
and define the seminorm $\parallel $ $\parallel $ on $A_\infty$
by $\parallel\phi_{k, \infty} (a)\parallel  = \lim \sup_l \parallel
 \varphi_\ell \circ
\dots \circ \phi_k (a) \parallel $ for $a \in A_k$.
Then the quotient of $A_\infty$ by the subspace  of elements with zero
seminorm becomes a normed space. The direct limit, $A$, of the system is
defined to be the completion of this normed space.
In fact the Banach space is matricially normed in the obvious way
and the maps $\phi_{k,\infty}$ are completely contractive.

If the $A_i$'s are
algebras then $A_\infty$ has a natural algebra structure,
the induced (operator) norm is an
algebra norm, and  the direct limit is a Banach
algebra.
In this case (see \cite{scp-simplicial}) $A$ is completely isomorphic to a
Hilbert space operator algebra.
If the $A_i$'s are C*-algebras and the $\varphi_i$'s are
contractive star homomorphisms (not necessarily injective), then $A_\infty$
inherits a natural star algebra structure, the norm defined above
is a C*-norm, and  the direct limit is a C*-algebra.
\vspace{.3in}

{\bf Lemma 9.}
\it \ The system above
can be replaced by a subsystem

\begin{center}
\setlength{\unitlength}{0.0125in}%
\begin{picture}(225,22)(125,755)
\thicklines
\put(240,760){\vector( 1, 0){ 40}}
\put(210,755){\makebox(0,0)[lb]{\raisebox{0pt}[0pt][0pt]{$\ca_2$}}}
\put(125,755){\makebox(0,0)[lb]{\raisebox{0pt}[0pt][0pt]{$\ca_1$}}}
\put(150,760){\vector( 1, 0){ 40}}
\put(310,755){\makebox(0,0)[lb]{\raisebox{0pt}[0pt][0pt]{\twlit . . . . }}}
\put(250,765){\makebox(0,0)[lb]{\raisebox{0pt}[0pt][0pt]{$\alpha_2$}}}
\put(160,765){\makebox(0,0)[lb]{\raisebox{0pt}[0pt][0pt]{$\alpha_1$}}}
\end{picture}
\end{center}

such that

(i) each $\alpha_i$ is of compression type,

(ii) the set of
${\cal A}_k$-irreducible projections associated with $\alpha_l
\circ \dots \circ
\alpha_k$ for each $ \ell > k$ is  the
same as the set of ${\cal A}_k$-irreducible
projections associated with $\alpha_k$,

(iii) the restriction of each map $\alpha_{k+1}$ to
$\alpha_k ({\cal A}_k)$ is isometric.
\rm
\vspace{.3in}

\begin{mproof}
A composition $\beta \circ \alpha$ of compression type maps is of compression
type.
Furthermore, the compression projections for the compositions are
subprojections
 of the
compression projections for $\alpha$ (by Proposition 8). In view of this we may
choose $k$ so that the compression projections for the composition
\[
\varphi_\ell \circ \dots \circ \varphi_1
\]
is constant for all $l > k$. It is now clear how to similarly choose
the maps $\alpha_k$ to satisfy the conditions (i) and (ii).
Property (iii) now follows.
\end{mproof}
\vspace{.3in}

{\bf Theorem  10.} \it \ Let
$\{A_k,\phi_k\}$ be a contractive regular direct system of digraph spaces
and let $\{\ca_k,\alpha_k\}$ be an essentially isometric subsystem
with the properties of Lemma 9.
Then there exists isometric star homomorphisms $\psi_k$ such that the diagrams

\begin{center}
\setlength{\unitlength}{0.0125in}%
\begin{picture}(155,107)(80,675)
\thicklines
\put(150,680){\vector( 1, 0){ 40}}
\put(115,740){\vector( 0,-1){ 40}}
\put(235,740){\vector( 0,-1){ 40}}
\put(220,760){\makebox(0,0)[lb]{\raisebox{0pt}[0pt][0pt]{$
 \alpha_{k+1}(\ca_{k+1})$}}}
\put(100,760){\makebox(0,0)[lb]{\raisebox{0pt}[0pt][0pt]{$
 \alpha_{k}(\ca_{k}$)}}}
\put(150,760){\vector( 1, 0){ 40}}
\put( 90,720){\makebox(0,0)[lb]{\raisebox{0pt}[0pt][0pt]{$ i_k$}}}
\put(200,675){\makebox(0,0)[lb]{\raisebox{0pt}[0pt][0pt]{$
 C^*(\alpha_{k+1}(\ca_{k+1}))$}}}
\put(210,720){\makebox(0,0)[lb]{\raisebox{0pt}[0pt][0pt]{$ i_{k+1}$}}}
\put(160,770){\makebox(0,0)[lb]{\raisebox{0pt}[0pt][0pt]{$ \alpha_{k+1}$}}}
\put(160,690){\makebox(0,0)[lb]{\raisebox{0pt}[0pt][0pt]{$ \psi_k$}}}
\put( 80,675){\makebox(0,0)[lb]{\raisebox{0pt}[0pt][0pt]{$
 C^*(\alpha_{k}(\ca_{k}))$}}}
\end{picture}
\end{center}

commute.
\rm
\vspace{.3in}

\begin{mproof}
By Theorem 5 we may assume that the maps are of compression type.
Let ${\cal A}_k \subseteq M_{m_k}$ be the natural inclusion.
Let $\alpha_k = \gamma = \sum_{r=1}^t \op \gamma_r $
where each $\gamma_r$ is an elementary compression type map associated with
diagonal projections $(Q_r^k, P_r^{k+1})$, and let $\alpha_{k+1} = \sum \op
\eta_s$ where $\eta_s$ is associated with diagonal projections
$(Q_s^{k+1}, P_s^{k+2})$. Let $I$ be a subset of indices such that
$\{Q_i^{k} : i \in I \}$ contains exactly one copy of each
$Q_r^k$ that appears in the description of $\alpha_k$.
Let $J$ be a subset of indices such that $\{Q_j^{k+1} : j \in J \}$
contains exactly one copy of each $Q_s^{k+1}$ that appears in the
description of $\alpha_{k+1}$.

By Proposition 7, we have the isomorphisms
\[
C^* (\alpha_k ({\cal A}_k)) \cong \sum_{i \in I} \oplus B (Q_i^k
\ \IC^{m_k})
\]
and
\[
C^* (\alpha_{k+1} ({\cal A}_{k+1})) \cong \sum_{j \in J} \oplus B
(Q_j^{k+1} \ \IC^{m_{k+1}}) \cong \sum_{j \in J} \oplus B (P_j^{k+2}
\ \IC^{m_{k+2}}).
\]

We shall define $\psi_k : C^* (\alpha_{k} (\ca_k)) \rightarrow C^*
(\alpha_{k+1} ({\cal A}_{k+1}))$ in accordance with the
multiplicity and the identification of
the embedding of each $Q_i^k$-summand into each
$P_j^{k+2}$-summand.

By Proposition 8, $\eta_j \circ \gamma = \sum \op \delta_\ell$ where
 $\delta_\ell$ is
an elementary compression type map associated with projections $(q_\ell
, p_\ell)$ where $\{p_\ell\}$ consists of mutually orthogonal
subprojections of $P_j^{k+2}$. By the hypotheses (that of Lemma 9 (ii)
in this case), each $q_\ell$ belongs to $ \{Q_i^k :
i \in I\}$. Note that if $Q_j^{k+1} = Q_m^{k+1}$ then the $\{q_\ell\}$
associated with $\eta_j \circ \gamma$ and the $\{q_\ell\}$ associated with
$\eta_m \circ \gamma$ are identical (counting multiple copies).

We are now ready to define $\psi_k.$ For $i \in I$ and $j \in J$, let
$n_{ij}$ be the number of copies of $Q_i^k$ that occur in the set of
$\{q_\ell\}$ in the description above of $\eta_j \circ  \gamma$.
Then there is a natural embedding of multiplicity $n_{ij}$ of the
$Q_i^k$-summand of
$C^* (\alpha_k (\ca_k))$ into the $P_j^{k+2}$-summand of
$C^* (\alpha_{k+1}(\ca_{k+1}))$.
This map is just the star extension of the restriction of
$\eta_j \circ \gamma$.

By assumption on the maps $\alpha_k$, given $i \in I$, there
exists $j \in J$ such that $n_{ij} \neq 0$. (Otherwise $\alpha_{k+1} \circ
 \gamma_i
= 0$ implying that $Q_i^k$ does not appear in the set of projections
associated with $\alpha_{k+1} \circ \alpha_k)$. Thus $\psi_k$ is isometric.

To see that the diagram commutes one can argue as follows. Consider an
element $b = \alpha_k(a)$ with $a \in \ca_k$. Then $b$ splits as a direct
sum $b_1 + \dots + b_t$ with $b_r = P_r^{k+1}b_rP_r^{k+1}$ where
$P_1^{k+1},\dots,P_t^{k+1}$ is the enumeration of all the projections in the
 definition of
$\alpha_k$. Furthermore, we can group the summands
according to the equivalence relation $ r \sim s$ on indices with
$Q_r^{k+1} = Q_s^{k+1}$ to obtain

\[
b = \sum_{i\in I} \op (\sum_{r\in i} \op b_r).
\]

Here we view the indices in $I$ also as the equivalence classes.
Each summand $b_r$ for $r \in i$ is the copy $\gamma_r(Q_r^kaQ_r^k)$
of $Q_i^kaQ_i^k$, and the map $i_k$ is the natural map which identifies these
copies. That is,
$i_k(b) = \sum_{i \in I} \op b_i$.

Similarly, each element $d$ of $\alpha_{k+1}(\ca_{k+1})$
has a direct sum representation

\[
d = \sum_{j\in J} \op (\sum_{s\in j} \op d_s).
\]

and $i_{k+1}$ is the map
with $i_{k+1}(d) = \sum_{j \in J} \op d_j$ which identifies multiple copies
as before.

Consider the representation of $\alpha_{k+1}(b)$ in the form
$\sum_{j\in J} \op (\sum_{s\in j} \op d_s)$.
Then each summand $d_s$, coming from $(\eta_s \circ \gamma)(b)$,
 is itself a direct sum of copies
of $b_i$ in the expression for $d_j$. The implication of this is that there
is a commuting diagram

\begin{center}
\setlength{\unitlength}{0.0125in}%
\begin{picture}(155,107)(80,675)
\thicklines
\put(150,680){\vector( 1, 0){ 40}}
\put(115,740){\vector( 0,-1){ 40}}
\put(235,740){\vector( 0,-1){ 40}}
\put(220,760){\makebox(0,0)[lb]{\raisebox{0pt}[0pt][0pt]{$
 \alpha_{k+1}(\ca_{k+1})$}}}
\put(100,760){\makebox(0,0)[lb]{\raisebox{0pt}[0pt][0pt]{$
 \alpha_{k}(\ca_{k}$)}}}
\put(150,760){\vector( 1, 0){ 40}}
\put( 90,720){\makebox(0,0)[lb]{\raisebox{0pt}[0pt][0pt]{$ i_k$}}}
\put(200,675){\makebox(0,0)[lb]{\raisebox{0pt}[0pt][0pt]
{$\{\sum_{j\in J}\op d_j\}$}}}
\put(210,720){\makebox(0,0)[lb]{\raisebox{0pt}[0pt][0pt]{$ i_{k+1}$}}}
\put(160,770){\makebox(0,0)[lb]{\raisebox{0pt}[0pt][0pt]{$ \alpha_{k+1}$}}}
\put(160,690){\makebox(0,0)[lb]{\raisebox{0pt}[0pt][0pt]{$ \hat{\psi}_{k+1}$}}}
\put( 80,675){\makebox(0,0)[lb]{\raisebox{0pt}[0pt][0pt]
{$\{\sum_{i\in I}\op b_i\}$}}}
\end{picture}
\end{center}

The map $\psi_{k+1}$ is simply the star-extension of the map
$ \hat{\psi}_{k+1}$ and so the diagram in the statement
of the theorem commutes.
\end{mproof}
\vspace{.3in}

{\bf Corollary 11.}\ \
\it
The direct limit of a contractive regular system of
digraph spaces
is completely isometric to a subspace of an AF
C*-algebra.
\rm
\vspace{.3in}

\begin{mproof}
In the notation above consider the following
commuting diagram for such a direct limit $A$.

\begin{center}
\setlength{\unitlength}{0.0125in}%
\begin{picture}(320,107)(80,675)
\thicklines
\put(280,680){\vector( 1, 0){ 40}}
\put(400,740){\vector( 0,-1){ 40}}
\put(150,760){\vector( 1, 0){ 40}}
\put(150,680){\vector( 1, 0){ 40}}
\put(115,740){\vector( 0,-1){ 40}}
\put(235,740){\vector( 0,-1){ 40}}
\put(220,760){\makebox(0,0)[lb]{\raisebox{0pt}[0pt][0pt]{$ \alpha_2(\ca_2)$}}}
\put(100,760){\makebox(0,0)[lb]{\raisebox{0pt}[0pt][0pt]{$ \alpha_1(\ca_1)$}}}
\put( 90,720){\makebox(0,0)[lb]{\raisebox{0pt}[0pt][0pt]{$ i_1$}}}
\put(280,760){\vector( 1, 0){ 40}}
\put(210,720){\makebox(0,0)[lb]{\raisebox{0pt}[0pt][0pt]{$ i_2$}}}
\put(200,675){\makebox(0,0)[lb]{\raisebox{0pt}[0pt][0pt]{$
 C^*(\alpha_2(\ca_2))$}}}
\put(340,760){\makebox(0,0)[lb]{\raisebox{0pt}[0pt][0pt]{$ . . . . $}}}
\put(340,720){\makebox(0,0)[lb]{\raisebox{0pt}[0pt][0pt]{$ . . . . $}}}
\put(340,680){\makebox(0,0)[lb]{\raisebox{0pt}[0pt][0pt]{$ . . . .$}}}
\put(400,760){\makebox(0,0)[lb]{\raisebox{0pt}[0pt][0pt]{$ A$}}}
\put(400,675){\makebox(0,0)[lb]{\raisebox{0pt}[0pt][0pt]{$ B$}}}
\put(385,720){\makebox(0,0)[lb]{\raisebox{0pt}[0pt][0pt]{$ i$}}}
\put(160,770){\makebox(0,0)[lb]{\raisebox{0pt}[0pt][0pt]{$ \alpha_2$}}}
\put(160,690){\makebox(0,0)[lb]{\raisebox{0pt}[0pt][0pt]{$ \psi_2$}}}
\put( 75,675){\makebox(0,0)[lb]{\raisebox{0pt}[0pt][0pt]{$
 C^*(\alpha_1(\ca_1))$}}}
\end{picture}
\end{center}

Since the inclusion maps $i_k$  are completely
isometric
we conclude that the induced
map $i: A \rightarrow B$ is a complete isometry.
\end{mproof}
\vspace{.3in}

{\bf Proposition 12.}\ \it Given the systems

\begin{center}
\setlength{\unitlength}{0.0125in}%
\begin{picture}(320,107)(80,675)
\thicklines
\put(280,680){\vector( 1, 0){ 40}}
\put(400,740){\vector( 0,-1){ 40}}
\put(150,760){\vector( 1, 0){ 40}}
\put(150,680){\vector( 1, 0){ 40}}
\put(115,740){\vector( 0,-1){ 40}}
\put(235,740){\vector( 0,-1){ 40}}
\put(220,760){\makebox(0,0)[lb]{\raisebox{0pt}[0pt][0pt]{$ \alpha_2(A_2)$}}}
\put(100,760){\makebox(0,0)[lb]{\raisebox{0pt}[0pt][0pt]{$ \alpha_1(A_1)$}}}
\put( 90,720){\makebox(0,0)[lb]{\raisebox{0pt}[0pt][0pt]{$ i_1$}}}
\put(280,760){\vector( 1, 0){ 40}}
\put(210,720){\makebox(0,0)[lb]{\raisebox{0pt}[0pt][0pt]{$ i_2$}}}
\put(200,675){\makebox(0,0)[lb]{\raisebox{0pt}[0pt][0pt]{$
 C^*(\alpha_2(A_2))$}}}
\put(340,760){\makebox(0,0)[lb]{\raisebox{0pt}[0pt][0pt]{$ . . . . $}}}
\put(340,720){\makebox(0,0)[lb]{\raisebox{0pt}[0pt][0pt]{$ . . . . $}}}
\put(340,680){\makebox(0,0)[lb]{\raisebox{0pt}[0pt][0pt]{$ . . . .$}}}
\put(400,760){\makebox(0,0)[lb]{\raisebox{0pt}[0pt][0pt]{$ A$}}}
\put(400,675){\makebox(0,0)[lb]{\raisebox{0pt}[0pt][0pt]{$ B$}}}
\put(385,720){\makebox(0,0)[lb]{\raisebox{0pt}[0pt][0pt]{$ i$}}}
\put(160,770){\makebox(0,0)[lb]{\raisebox{0pt}[0pt][0pt]{$ \alpha_2$}}}
\put(160,690){\makebox(0,0)[lb]{\raisebox{0pt}[0pt][0pt]{$ \psi_2$}}}
\put( 75,675){\makebox(0,0)[lb]{\raisebox{0pt}[0pt][0pt]{$
 C^*(\alpha_1(A_1))$}}}
\end{picture}
\end{center}

as in Theorem 10, define $f_k = \alpha_{k, \infty}
(I_k)$ where $I_k$ is the identity in ${\cal A}_k$. Then

(a) the sequence $e_k = i(f_k)$ forms a norm approximate identity for $B$

(b) the following statements are equivalent
 \begin{itemize}
\item[i)]  $\psi_i$ is unital for all large $i$.
 \item[ii)] $B$ has a unit $e$.
\item[iii)] $e_k$ converges in norm to an element $ e$ in  $B$.
\item[iv)] $f_k$ converges in norm to an element $g$ in $A$.
\end{itemize}
Under any of the assumptions in (b) we have that $i(g) = e$.
\rm
\vspace{.3in}

\begin{mproof}
The statements follow directly once it is observed that $i_k
(\alpha_k (I_k))$ is the identity in $B_k$ and that $i \circ \alpha_{k+1,
\infty} = \psi_{k, \infty} \circ i_k$.
\end{mproof}
\vspace{.3in}

{\bf Definition 13.} Let

\begin{center}
\setlength{\unitlength}{0.0125in}%
\begin{picture}(225,22)(125,755)
\thicklines
\put(240,760){\vector( 1, 0){ 40}}
\put(210,755){\makebox(0,0)[lb]{\raisebox{0pt}[0pt][0pt]{$A_2$}}}
\put(125,755){\makebox(0,0)[lb]{\raisebox{0pt}[0pt][0pt]{$A_1$}}}
\put(150,760){\vector( 1, 0){ 40}}
\put(310,755){\makebox(0,0)[lb]{\raisebox{0pt}[0pt][0pt]{\twlit . . . . }}}
\put(250,765){\makebox(0,0)[lb]{\raisebox{0pt}[0pt][0pt]{$\phi_2$}}}
\put(350,755){\makebox(0,0)[lb]{\raisebox{0pt}[0pt][0pt]{$ A$}}}
\put(160,765){\makebox(0,0)[lb]{\raisebox{0pt}[0pt][0pt]{$\phi_1$}}}
\end{picture}
\end{center}

 be a direct system of digraph spaces
where each $\varphi_i$ is a contractive regular bimodule
map. This system is said to be
{\em essentially unital} if the telescoped system
given by Lemma 9 satisfies any of the equivalent properties in
(b) of Proposition 12.
\vspace{.3in}

Before going further we provide a brief review of the idea
of a C*-envelope.

The appropriate setting, considered by Hamana (\cite{ham}),
is the category of unital
operator spaces (that is, self-adjoint unital subspaces of C*-algebras),
together with unital complete order injections. If $A$ is an
operator space, then a C*-{\it extension} of $A$ is a C*-algebra $B$
together with a unital complete order injection $\rho$ of $A$ into $B$
 such that $C^*(\rho(A)) = B$. A C*-extension $B$ is a C*-{\it envelope}
 of $A$ provided that, given any operator system, $C$, and any unital
 completely positive map $\tau: B \rightarrow C$, $\tau$ is a complete
 order injection whenever $\tau \circ \rho$ is. Hamana
 proves the existence and uniqueness (up to a suitable notion of
 equivalence) of C*-envelopes. Furthermore , he shows that the C*-envelope
 of $A$ is a minimal C*-extension in the family of all C*-extensions
 of $A$.

After proving the existence of C*-envelopes, Hamana then uses this
to prove the existence of a Silov boundary for $A$,
a theme first developed by Arveson  \cite{arv}.
The Silov boundary is a generalization to operator spaces,
of the
usual notion of Silov boundary from function spaces.

Let $B$ be a
C*-algebra and let $A$ be a unital sub-operator-space such that
$B = C^*(A)$. An ideal $J$ in $B$ is called a {\it boundary ideal} for
$A$ if the canonical quotient map $B \rightarrow B/J$ is completely
isometric on $A$. A boundary ideal exists which contains every other boundary
ideal and this maximal boundary ideal
is called the {\it Silov boundary} for $A$. Hamana shows that
if $B$ is a C*-extension for $A$, then the C*-envelope for $A$ is
isomorphic to $B/J$, where $J$ is the Silov boundary for $A$.
It is this form of the definition of the
C*-envelope that will be used below.

If $A$ is merely a unital subspace of a C*-algebra, rather than
an operator space, then we define the C*-envelope of $A$ to be the C*-envelope
of the operator space $A + A^*$.
\vspace{.3in}

{\bf Proposition 14.}\
\it Let $A$ be the direct limit of an essentially unital
contractive regular
direct system of digraph spaces.
Then $A$ is completely isometrically isomorphic
to a unital subspace of an
AF C*-algebra. Furthermore the C*-envelope of $A$
is an AF C*-algebra.
\vspace{.3in}
\rm

\begin{mproof}
Consider the telescoped system and the associated commuting diagram

\begin{center}
\setlength{\unitlength}{0.0125in}%
\begin{picture}(320,107)(80,675)
\thicklines
\put(280,680){\vector( 1, 0){ 40}}
\put(400,740){\vector( 0,-1){ 40}}
\put(150,760){\vector( 1, 0){ 40}}
\put(150,680){\vector( 1, 0){ 40}}
\put(115,740){\vector( 0,-1){ 40}}
\put(235,740){\vector( 0,-1){ 40}}
\put(220,760){\makebox(0,0)[lb]{\raisebox{0pt}[0pt][0pt]{$ \alpha_2(\ca_2)$}}}
\put(100,760){\makebox(0,0)[lb]{\raisebox{0pt}[0pt][0pt]{$ \alpha_1(\ca_1)$}}}
\put( 90,720){\makebox(0,0)[lb]{\raisebox{0pt}[0pt][0pt]{$ i_1$}}}
\put(280,760){\vector( 1, 0){ 40}}
\put(210,720){\makebox(0,0)[lb]{\raisebox{0pt}[0pt][0pt]{$ i_2$}}}
\put(225,675){\makebox(0,0)[lb]{\raisebox{0pt}[0pt][0pt]{$ B_2$}}}
\put(340,760){\makebox(0,0)[lb]{\raisebox{0pt}[0pt][0pt]{$ . . . . $}}}
\put(340,720){\makebox(0,0)[lb]{\raisebox{0pt}[0pt][0pt]{$ . . . . $}}}
\put(340,680){\makebox(0,0)[lb]{\raisebox{0pt}[0pt][0pt]{$ . . . .$}}}
\put(400,760){\makebox(0,0)[lb]{\raisebox{0pt}[0pt][0pt]{$ A$}}}
\put(400,675){\makebox(0,0)[lb]{\raisebox{0pt}[0pt][0pt]{$ B$}}}
\put(385,720){\makebox(0,0)[lb]{\raisebox{0pt}[0pt][0pt]{$ i$}}}
\put(160,770){\makebox(0,0)[lb]{\raisebox{0pt}[0pt][0pt]{$ \alpha_2$}}}
\put(160,690){\makebox(0,0)[lb]{\raisebox{0pt}[0pt][0pt]{$ \psi_2$}}}
\put(105,675){\makebox(0,0)[lb]{\raisebox{0pt}[0pt][0pt]{$ B_1$}}}
\end{picture}
\end{center}

where $B_i = C^* (\alpha_i (\ca_i))$
is a finite-dimensional C*-algebra, each
$\psi_i$ is an isometric star homomorophism and $i$ is a complete
isometry.
Since the telescoped system is unital,
we see from the definition of $\psi_i$ that
 each $\psi_i$ is unital. Thus, by Proposition 12
 there exists  $g \in A$ such
that $i(g) = e$, the unit in $B$ and the first assertion of the theorem
follows.

 Let $A$ also denote the image of $A$
in $B$. Then $\tilde{A} = \overline{A+A^*}$ is an operator space and,
since $i$ is a complete isometry on $A$ it extends to
a unital complete order injection $\ i: \tilde{A} \to B$. Thus
$B$ is a C*-extension of
$\tilde{A}$  and $ C^*_{env} (A)$ is (completely isometric to) the quotient
$B/J$
where $J$ is the Silov boundary of $\tilde{A}$ in $B$.
\end{mproof}
\vspace{.3in}

\bf The direct system for $ C^*_{env} (A)$ \rm
\vspace{.3in}

We shall now identify the direct system for $ C^*_{env} (A)$. This requires an
identification of the Silov ideal  $J$ above, which, as we see in
Example 16 below,
may be nonzero. This example shows
that we cannot expect $C^*_{env} (A)$ to be an isometric limit
of the C*-algebras $C^*(\alpha_k(A_k))$ of the isometric
telescoped system for $A$. Nevertheless, because of the
nature of ideals in direct systems it follows that the
C*-envelope is a direct limit of these building blocks
with respect to not necessarily injective embeddings.

View $B_k$ as the isometric
image of $B_k$ in $B$ and let $J_k = J \bigcap B_k$. Then we have $J =
\overline{\bigcup J_k}$ and isometric *-homomorphisms
$\tilde{\psi_i}$ exist
such that the following diagram commutes. (See Bratteli \cite{bra}
for example.)

\begin{center}
\setlength{\unitlength}{0.0125in}%
\begin{picture}(320,107)(80,675)
\thicklines
\put(280,680){\vector( 1, 0){ 40}}
\put(400,740){\vector( 0,-1){ 40}}
\put(150,760){\vector( 1, 0){ 40}}
\put(150,680){\vector( 1, 0){ 40}}
\put(115,740){\vector( 0,-1){ 40}}
\put(235,740){\vector( 0,-1){ 40}}
\put(220,760){\makebox(0,0)[lb]{\raisebox{0pt}[0pt][0pt]{$ B_2$}}}
\put(100,760){\makebox(0,0)[lb]{\raisebox{0pt}[0pt][0pt]{$ B_1$}}}
\put( 90,720){\makebox(0,0)[lb]{\raisebox{0pt}[0pt][0pt]{$ \pi_1$}}}
\put(280,760){\vector( 1, 0){ 40}}
\put(210,720){\makebox(0,0)[lb]{\raisebox{0pt}[0pt][0pt]{$ \pi_2$}}}
\put(225,675){\makebox(0,0)[lb]{\raisebox{0pt}[0pt][0pt]{$ B_2/J_2$}}}
\put(340,760){\makebox(0,0)[lb]{\raisebox{0pt}[0pt][0pt]{$ . . . . $}}}
\put(340,720){\makebox(0,0)[lb]{\raisebox{0pt}[0pt][0pt]{$ . . . . $}}}
\put(340,680){\makebox(0,0)[lb]{\raisebox{0pt}[0pt][0pt]{$ . . . .$}}}
\put(400,760){\makebox(0,0)[lb]{\raisebox{0pt}[0pt][0pt]{$ B$}}}
\put(400,675){\makebox(0,0)[lb]{\raisebox{0pt}[0pt][0pt]{$ B/J$}}}
\put(385,720){\makebox(0,0)[lb]{\raisebox{0pt}[0pt][0pt]{$ \pi$}}}
\put(160,770){\makebox(0,0)[lb]{\raisebox{0pt}[0pt][0pt]{$ \psi_1$}}}
\put(160,690){\makebox(0,0)[lb]{\raisebox{0pt}[0pt][0pt]{$ \tilde{\psi}_1$}}}
\put(105,675){\makebox(0,0)[lb]{\raisebox{0pt}[0pt][0pt]{$ B_1/J_1$}}}
\end{picture}
\end{center}

Since $J$ is a boundary ideal of $\tilde{A}$ in $B$,
the map $\pi : B \rightarrow
B/J$ is completely isometric on $\tilde{A}$. Hence $\pi_k : B_k
\rightarrow B_k / J_k$ is completely isometric on $i_k (\alpha_k
(A_k))$. We thus have the commuting system

\begin{center}
\setlength{\unitlength}{0.0125in}%
\begin{picture}(320,107)(80,675)
\thicklines
\put(280,680){\vector( 1, 0){ 40}}
\put(400,740){\vector( 0,-1){ 40}}
\put(150,760){\vector( 1, 0){ 40}}
\put(150,680){\vector( 1, 0){ 40}}
\put(115,740){\vector( 0,-1){ 40}}
\put(235,740){\vector( 0,-1){ 40}}
\put(220,760){\makebox(0,0)[lb]{\raisebox{0pt}[0pt][0pt]{$ \alpha_2(A_2)$}}}
\put(100,760){\makebox(0,0)[lb]{\raisebox{0pt}[0pt][0pt]{$ \alpha_1(A_1)$}}}
\put( 90,720){\makebox(0,0)[lb]{\raisebox{0pt}[0pt][0pt]{$ \tilde{i}_1$}}}
\put(280,760){\vector( 1, 0){ 40}}
\put(210,720){\makebox(0,0)[lb]{\raisebox{0pt}[0pt][0pt]{$ \tilde{i}_2$}}}
\put(200,675){\makebox(0,0)[lb]{\raisebox{0pt}[0pt][0pt]{$ B_2/J_2$}}}
\put(340,760){\makebox(0,0)[lb]{\raisebox{0pt}[0pt][0pt]{$ . . . . $}}}
\put(340,720){\makebox(0,0)[lb]{\raisebox{0pt}[0pt][0pt]{$ . . . . $}}}
\put(340,680){\makebox(0,0)[lb]{\raisebox{0pt}[0pt][0pt]{$ . . . .$}}}
\put(400,760){\makebox(0,0)[lb]{\raisebox{0pt}[0pt][0pt]{$ A$}}}
\put(400,675){\makebox(0,0)[lb]{\raisebox{0pt}[0pt][0pt]{$ B/J =
 C^*_{env}(A)$}}}
\put(385,720){\makebox(0,0)[lb]{\raisebox{0pt}[0pt][0pt]{$ \tilde{i}$}}}
\put(160,770){\makebox(0,0)[lb]{\raisebox{0pt}[0pt][0pt]{$ \alpha_2$}}}
\put(160,690){\makebox(0,0)[lb]{\raisebox{0pt}[0pt][0pt]{$ \tilde{\psi}_2$}}}
\put( 75,675){\makebox(0,0)[lb]{\raisebox{0pt}[0pt][0pt]{$ B_1/J_1$}}}
\end{picture}
\end{center}

where $\tilde{i_k} = \pi_k \circ i_k$ is a complete isometry and
$\tilde{i} = \pi \circ i$ is a complete isometry.

Furthermore, as Ken Davidson has noted,
the identifications of Theorem 10 allow
the Silov boundary $J$ to be identified in the
following specific intrinsic manner.

Consider the summands of $C^*(\alpha_k(A_k))$ which are {\em not maximal}
in the sense that they correspond to proper compressions of larger
summands of $C^*(\alpha_k(A_k))$. Otherwise refer to a summand as {\em
maximal}.
If a summand of $C^*(\alpha_k(A_k))$ is never mapped into a maximal summand
of $C^*(\alpha_j(A_j))$ for any $j > k$ then (the image of)
this summand is clearly
contained in the Silov ideal. (In fact such a summand generates
a boundary ideal.) Moreover, the Silov boundary $J$ is precisely the ideal, $K$
say, that is generated by all such summands.

To see this suppose, by way of contradiction, that $J$ contains $K$ strictly.
Then there is a summand, $S$ say,
of $C^*(\alpha_k(A_k))$ for some $k$, which is contained in $J$
and which has (partial)
embeddings into maximal summands
$M_{n_j} \subseteq C^*(\alpha_{n_j}(A_{n_j}))$, for
some increasing sequence $n_j$.  Let $a \in \alpha_k(A_k)$
be an element such that $\|i_k(a) + S\| < \|a\|$. For example,
pick  $a = \alpha_k(b)$ where $b$ is
supported by a  single compression projection
(for the domain of $\alpha_k$) corresponding to $S$ and all
proper subcompressions
of $b$ have strictly smaller norm.
Plainly $\|\ga_{k+1,\infty}(\tilde{i}(a)) + J\|
\le \|i_k(a) + S\| <  \|a\|$. But the existence of embeddings
into maximal summands implies that
$\|\ga_{k+1,\infty}(\tilde{i}(a))\| = \|a\|$. This is contrary to the fact that
$J$ is a boundary ideal and so the assertion follows.

In summary then we have the following theorem.
\vspace{.3in}

{\bf Theorem 15.} \ Let $A$ be the direct limit of an essentially unital
contractive regular direct
system of digraph spaces $A_k$. Then
\begin{enumerate}
\item[(i)] $A$ is completely isometrically isomorphic to a unital
subspace of an AF $C^*$-algebra.
\item[(ii)] $C^*_{env}(A)$ is an AF  $C^*$-algebra.
\item[(iii)] If $B = \displaystyle{\lim_\to} C^*(\ga_i(\ca_i)) $
is the  $C^*$-algebra of the telescoped system for $A$ (as in Theorem 10) then
$C^*_{env}(A) = B/J$ where $J$ is the ideal generated by those
summands of $C^*(\ga_i(\ca_i))$, for $i = 1,2,\dots, $ which have no
partial embedding into a maximal summand of $C^*(\ga_j(\ca_j))$
for all $j > i$.
\end{enumerate}

{\bf Remark.} \ Note that if $\varphi : A (G) \rightarrow A(H)$ is a
contractive
regular bimodule map, then since $\varphi$ is of compression type,
$\varphi$ can be extended to $\tilde{\varphi} : A(G) + A(G)^* \rightarrow
A(H) + A(H)^*$ where $\tilde{\varphi}$ is of compression type associated
with the same projections as $\varphi$.
It follows that the
commuting diagram above
can be interpolated to yield

\begin{center}
\setlength{\unitlength}{0.0125in}%
\begin{picture}(395,149)(30,595)
\thicklines
\put(360,660){\vector( 1, 0){ 40}}
\put(360,600){\vector( 1, 0){ 40}}
\put(265,720){\vector( 1, 0){ 40}}
\put(265,660){\vector( 1, 0){ 40}}
\put(265,600){\vector( 1, 0){ 40}}
\put(120,720){\vector( 1, 0){ 40}}
\put(120,660){\vector( 1, 0){ 40}}
\put(120,600){\vector( 1, 0){ 40}}
\put( 95,700){\vector( 0,-1){ 20}}
\put( 95,640){\vector( 0,-1){ 20}}
\put(200,705){\vector( 0,-1){ 20}}
\put(200,640){\vector( 0,-1){ 20}}
\put(425,700){\vector( 0,-1){ 20}}
\put(360,720){\vector( 1, 0){ 40}}
\put(425,640){\vector( 0,-1){ 20}}
\put(175,660){\makebox(0,0)[lb]{\raisebox{0pt}[0pt][0pt]{$
\ga_2(A_2+A_2^*)$}}}
\put(185,715){\makebox(0,0)[lb]{\raisebox{0pt}[0pt][0pt]{$  \ga_2(A_2)$}}}
\put(130,735){\makebox(0,0)[lb]{\raisebox{0pt}[0pt][0pt]{$  \ga_2$}}}
\put(130,670){\makebox(0,0)[lb]{\raisebox{0pt}[0pt][0pt]{$  \tilde{\ga}_2$}}}
\put(130,610){\makebox(0,0)[lb]{\raisebox{0pt}[0pt][0pt]{$  \tilde{\psi_1}$}}}
\put( 70,715){\makebox(0,0)[lb]{\raisebox{0pt}[0pt][0pt]{$  \ga_1(A_1)$}}}
\put(315,660){\makebox(0,0)[lb]{\raisebox{0pt}[0pt][0pt]{  .......}}}
\put(425,715){\makebox(0,0)[lb]{\raisebox{0pt}[0pt][0pt]{$  A$}}}
\put(420,655){\makebox(0,0)[lb]{\raisebox{0pt}[0pt][0pt]{$  \tilde{A}$}}}
\put(420,595){\makebox(0,0)[lb]{\raisebox{0pt}[0pt][0pt]{$  B/J =
 C^*_{env}(A)$}}}
\put( 45,660){\makebox(0,0)[lb]{\raisebox{0pt}[0pt][0pt]{$
\ga_1(A_1+A_1^*)$}}}
\put( 30,595){\makebox(0,0)[lb]{\raisebox{0pt}[0pt][0pt]{$
 C^*(\ga_1(A_1))/J_1$}}}
\put(170,595){\makebox(0,0)[lb]{\raisebox{0pt}[0pt][0pt]{$
 C^*(\ga_2(A_2))/J_2$}}}
\end{picture}
\end{center}

Consequently the operator space
$\tilde{A} \equiv \overline{i (A) + i (A)^*}$ is completely
isometric to the limit of the middle system.
\vspace{.3in}

{\bf Example 16} \ \ The following example illustates
how the ideal $J$ of Theorem 15 (iii) may be proper,
even when the system $\{C^*(\ga_i(\ca_i)), \psi_i\}$ is isometric.

Consider
\begin{center}
\setlength{\unitlength}{0.0125in}%
\begin{picture}(225,22)(125,755)
\thicklines
\put(240,760){\vector( 1, 0){ 40}}
\put(210,755){\makebox(0,0)[lb]{\raisebox{0pt}[0pt][0pt]{$T_5$}}}
\put(125,755){\makebox(0,0)[lb]{\raisebox{0pt}[0pt][0pt]{$T_3$}}}
\put(150,760){\vector( 1, 0){ 40}}
\put(310,755){\makebox(0,0)[lb]{\raisebox{0pt}[0pt][0pt]{\twlit . . . . }}}
\put(250,765){\makebox(0,0)[lb]{\raisebox{0pt}[0pt][0pt]{$\alpha_2$}}}
\put(350,755){\makebox(0,0)[lb]{\raisebox{0pt}[0pt][0pt]{$ A$}}}
\put(160,765){\makebox(0,0)[lb]{\raisebox{0pt}[0pt][0pt]{$\alpha_1$}}}
\end{picture}
\end{center}
where $T_n$ is the
algebra of $n \times n$
upper triangular matrices and where $\alpha_i (a) = a_{22} \oplus pap
\oplus pap$ where $a_{22}$ is the compression of $a$ to the second
minimal diagonal projection and
$p$ is compression to the last $n-1$ minimal diagonal
projections.

For example,
\[
\left[ \begin{array}{ccc}
a_{11} & a_{12} & a_{13} \\ & a_{22} & a_{23}\\
& & a_{33} \end{array} \right]
\rightarrow
\left[ \begin{array}{ccc}
a_{22} & &  \\ &C& \\ & & C \end{array}
\right]
= \alpha(T_3) \hspace{1ex} where \hspace{1ex} C = \left[
\begin{array}{cc} a_{22}&a_{23}\\ &a_{33} \end{array} \right],
\]

\[ \left[ \begin{array}{ccccc} a_{11}&a_{12}&a_{13}&a_{14}&a_{15}\\
&a_{22}&a_{23}&a_{24}&a_{25}\\& & a_{33}&a_{34}&a_{35}\\
& & & a_{44}&a_{45}\\& & & &a_{55} \end{array} \right] \rightarrow
\left[ \begin{array}{ccc} a_{22}& & \\ & D & \\ & & D \end{array}
\right] = \alpha_2(T_5) \hspace{1ex} \]
 where $\displaystyle D =
\left[ \begin{array}{cccc} a_{22}&a_{23}&a_{24}&a_{25}\\
& a_{33}&a_{34}&a_{35}\\& & a_{43} &a_{45}\\ & & & a_{55} \end{array}
\right],
$ and thus

\[
\left[ \begin{array}{ccc}
a_{11} & a_{12} & a_{13} \\ & a_{22} & a_{23}\\
& & a_{33} \end{array} \right]
\rightarrow
\left[ \begin{array}{ccccc} a_{22}& & & & \\ & C& & & \\ & & C & &\\
& & & C &\\ & & & & C \end{array} \right] = (\alpha_2 \circ \alpha_1)(T_3)
\]
 where $C$ is as above.

Note that the system satisfies the properties of Lemma 9. Let $Q_1$ be the
second minimal diagonal projection in $T_3$ and let $Q_2$ be the sum of
the last two minimal diagonal projections in $T_3$. Let $P_1$ be the first
minimal diagonal projection in $T_9$, $P_2$  the sum of the second
through the fifth minimal diagonal projections in $T_9$, and $P_3$
the sum of the sixth through the ninth minimal diagonal projections in
$T_9$. We can view $\alpha_2 \circ \alpha_1$ as embedding
$B(Q_1 \ \IC^3)$ into $B(P_1 \ \IC^9)$ with multiplicity 1 and embedding
$B(Q_2 \ \IC^3)$ into $B(P_2 \ \IC^9)$ with multiplicity 2 (as well as
$B(Q_2 \ \IC^3)$ into $B(P_3 \ \IC^9)$ with multiplicity 2).

Letting $M_i$ denote the $i \times i$ matrices, we have
\[C^*(\alpha_1(T_3)) \cong B(Q_1 \ \IC^3) \oplus B(Q_2 \ \IC^3)
\cong M_1 \oplus M_2\]
and
\[C^*(\alpha_2(T_5)) \cong B(P_1 \ \IC^9) \oplus B(P_2 \ \IC^9)
\cong M_1 \oplus M_4.\]
The embedding $\psi_2: C^*(\alpha_1(T_3)) \rightarrow C^*(\alpha_2(T_5))$
induced by $\alpha_2 \circ \alpha_1$ is represented by the Bratteli
diagram
\[ \begin{array}{ccc} M_1&\oplus&M_2\\ \mid& &\parallel \\
M_1&\oplus&M_4 \end{array}. \]
This Bratteli diagram format continues for
\begin{center}
\setlength{\unitlength}{0.0125in}%
\begin{picture}(225,22)(125,755)
\thicklines
\put(240,760){\vector( 1, 0){ 40}}
\put(210,755){\makebox(0,0)[lb]{\raisebox{0pt}[0pt][0pt]{$B_2$}}}
\put(125,755){\makebox(0,0)[lb]{\raisebox{0pt}[0pt][0pt]{$B_1$}}}
\put(150,760){\vector( 1, 0){ 40}}
\put(310,755){\makebox(0,0)[lb]{\raisebox{0pt}[0pt][0pt]{\twlit . . . . }}}
\put(250,765){\makebox(0,0)[lb]{\raisebox{0pt}[0pt][0pt]{$\psi_2$}}}
\put(350,755){\makebox(0,0)[lb]{\raisebox{0pt}[0pt][0pt]{$ B$}}}
\put(160,765){\makebox(0,0)[lb]{\raisebox{0pt}[0pt][0pt]{$\psi_1$}}}
\end{picture}
\end{center}
where $B_i \cong M_1 \oplus M_{2^i}$.

Let $J$ be the ideal in $B = \displaystyle{\lim_\to (B_i, \psi_i)}$
that corresponds to the subdiagram of embeddings of $M_1\oplus 0$
into $M_1 \oplus 0$.
  Then $\pi : B \rightarrow B/J$ is completely isometric on $A$
(deleting the corner entry
$a_{22}$ does not affect the norm of any of the images) and so $J$ is a
boundary ideal for $A$.
 and $B$ is not the C*-envelope of $A$.

In fact, as we see from the general discussion below, we can
make the identifications
$C^*_{env} (A) \cong B/J =  \displaystyle{\lim_\to (B_i/J_i,
\tilde{\psi_i)}}$
where $J_i = M_1 \oplus 0 \subseteq B_i$  and so $C^*_{env}
(A)$ is the UHF$(2^\infty)$ Glimm algebra.
\vspace{.3in}

{\bf Final remarks}
\vspace{.3in}

We conclude with comments on various contractive regular systems.

Consider the system

\begin{center}
\setlength{\unitlength}{0.0125in}%
\begin{picture}(225,22)(125,755)
\thicklines
\put(240,760){\vector( 1, 0){ 40}}
\put(200,755){\makebox(0,0)[lb]{\raisebox{0pt}[0pt][0pt]{$M_{n+k_1}$}}}
\put(125,755){\makebox(0,0)[lb]{\raisebox{0pt}[0pt][0pt]{$M_n$}}}
\put(150,760){\vector( 1, 0){ 40}}
\put(310,755){\makebox(0,0)[lb]{\raisebox{0pt}[0pt][0pt]{\twlit . . . . }}}
\put(250,765){\makebox(0,0)[lb]{\raisebox{0pt}[0pt][0pt]{$\phi_2$}}}
\put(350,755){\makebox(0,0)[lb]{\raisebox{0pt}[0pt][0pt]{$ B$}}}
\put(160,765){\makebox(0,0)[lb]{\raisebox{0pt}[0pt][0pt]{$\phi_1$}}}
\end{picture}
\end{center}

with embeddings $\phi_j$ such that
\[
\phi_j(a) = a \oplus a_{p_jp_j}I_{k_j}
\]
where $a_{p_jp_j}$ is the last diagonal entry of the matrix $a$.
Here $k_j$ is a sequence of positive integers and
$p_j = n + k_1 + \dots +  k_{j-1}.$  \ These embeddings
 restrict to algebra injections of the
upper triangular matrix subalgebras, giving a
triangular limit algebra,
the limit algebra  $A$ of Example B in \cite{hop-lau}.

In a sense this example is not properly of compression type
because the limit space $B$ can be viewed as a limit
of a subsystem for which the embeddings restrict to star
algebra homomorphisms.
Indeed, let $B_j \subseteq M_{p_j}$ be the block
diagonal  subspace $M_{p_j-1} \oplus \ \IC$. Then
the maps $\phi_j : B_{j-1} \to B_j$  are C*-algebra embeddings, and,
furthermore, the subsystem $\{B_{j-1},\phi_j\}$ has the same limit, $B$.
In particular the triangular operator algebra $A$ is a regular subalgebra
of an AF C*-algebra in the usual sense (\cite{scp-book}).
Here $B = C^*(A)$ is the C*-envelope of $A$ and,
furthermore, the masa $A \cap A^*$
in $A$ is a masa in $C^*_{env}(A)$. Similar remarks apply to
the Examples D, E, F of \cite{hop-lau}.

On the other hand (as noted in \cite{hop-lau}) the system
$ T_n \to T_{n+1} \to T_{n+2} \dots$ with algebra homomorphisms
$a \to a \oplus a_{ii}$, where $ 1 < i < n$ is fixed,
is properly of compression type; the containing system
$ M_n \to M_{n+1} \to M_{n+2} \dots$
does not have an algebra subsystem (in the sense above)
with the same limit. Indeed in this case the image of the
masa $A\cap A^*$ in $A$ under the inclusion $A \to C^*_{env}(A)$ is not
maximal abelian. For this example one readily identifies the C*-envelope
as the AF C*-algebra with Bratteli diagram

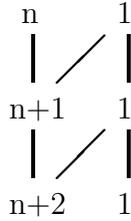
\begin{center}
\setlength{\unitlength}{0.0125in}%
\begin{picture}(50,89)(150,680)
\thicklines
\put(200,755){\line( 0,-1){ 20}}
\put(200,715){\line( 0,-1){ 20}}
\put(160,715){\line( 0,-1){ 20}}
\put(190,755){\line(-1,-1){ 20}}
\put(190,715){\line(-1,-1){ 20}}
\put(160,755){\line( 0,-1){ 20}}
\put(150,720){\makebox(0,0)[lb]{\raisebox{0pt}[0pt][0pt]{\twlrm n+1}}}
\put(195,680){\makebox(0,0)[lb]{\raisebox{0pt}[0pt][0pt]{\twlrm 1}}}
\put(150,680){\makebox(0,0)[lb]{\raisebox{0pt}[0pt][0pt]{\twlrm n+2}}}
\put(155,760){\makebox(0,0)[lb]{\raisebox{0pt}[0pt][0pt]{\twlrm n}}}
\put(195,760){\makebox(0,0)[lb]{\raisebox{0pt}[0pt][0pt]{\twlrm 1}}}
\put(195,720){\makebox(0,0)[lb]{\raisebox{0pt}[0pt][0pt]{\twlrm 1}}}
\end{picture}
\end{center}
\vspace{.3in}

In principle, the enveloping C*-algebra of a given (proper)
compression type system of digraph spaces can be identified by explicating
the Bratteli diagram from the process described in the Lemma 9, Theorem 14
and the related discussions. Nevertheless,  combinatorial
counting arguments may be necessary for this which make this
difficult in practice, as is the case, for example, with the system
$T_2 \to T_3 \to T_6 \to \dots$, with unital embedding homomorphisms
that have exactly one copy of every proper interval compression.

It should be apparent from our discussions that the
major aspect determining the C*-envelope is the nature
and the number of
the compressions appearing in the morphisms of the given direct system.
The identity of the building blocks themselves plays a minor role.
In fact, for any given AF C*-algebra $B$ one can construct a
regular contractive system $T_{2^{N_1}} \to T_{2^{N_2}} \to \dots$,
with algebra homomorphisms,
such that the triangular limit algebra has C*-envelope equal to $B$.

To see this let $\phi^\pr : M_{n_1} \op \dots \op M_{n_p} \to
M_{m_1} \op \dots \op M_{m_q} $ be any standard C*-algebra homomorphism and let
$\theta_1 : A_1 \to A_2$   be the restriction
of $\phi^\pr$ to the upper triangular
subalgebras. Choose $N_1, N_2$ with
$2^{N_1} \ge n_1 + \dots +n_p$,
$2^{N_2} \ge m_1 + \dots +m_q$ and
let $\kappa_i : T_{2^{N_i}} \to A_i$ be the natural
block diagonal compression maps. Then the map
$\ga_1 = \kappa_2^{-1} \circ \theta_1 \circ \kappa_1$
is a contractive regular algebra homomorphism from $A_1$ to $A_2$.
Iterating this construction one can
express $ \displaystyle{\lim_\to}(A_k,\theta_k)$
as $\displaystyle{\lim_\to}(T_{2^{N_k}},\ga_k)$ and the assertion follows.

\end{document}